\documentclass[leqno]{amsart}
\usepackage{amsfonts}
\usepackage{amscd,amsmath,amsopn,amssymb,amsthm,multicol}

\newtheorem{theorem}{Theorem}

\newtheorem{corollary}[theorem]{Corollary}

\newtheorem{definition}[theorem]{Definition}
\newtheorem{example}[theorem]{Example}

\newtheorem{lemma}[theorem]{Lemma}

\newtheorem{question}[theorem]{Question}
\newtheorem{proposition}[theorem]{Proposition}
\newtheorem{remark}[theorem]{Remark}

\input{tcilatex}

\begin{document}
\title{Covering $\mathbb{R}$-trees, $\mathbb{R}$-free groups, and Dendrites }
\author{V.N.~Berestovski\u\i\ and C.P.~Plaut}
\address{ Valeri\u\i\ Nikolaevich Berestovski\u\i\ \\
Omsk Branch of Sobolev Institute of Mathematics SD RAS \\
644099, Omsk, ul. Pevtsova, 13, Russia}
\email{berestov@ofim.oscsbras.ru}
\address{Conrad Peck Plaut \\
Department of Mathematics\\
University of Tennessee\\
Knoxville, TN 37919, USA}
\email{cplaut@math.utk.edu}

\begin{abstract}
We prove that every length space $X$ is the orbit space (with the quotient
metric) of an $\mathbb{R}$-tree $\overline{X}$ via a free action of a
locally free subgroup $\Gamma (X)$ of isometries of $\overline{X}$. The
mapping $\overline{\phi }:\overline{X}\rightarrow X$ is a kind of
generalized covering map called a URL-map and is universal among URL-maps
onto $X$. $\overline{X}$ is the unique $\mathbb{R}$-tree admitting a URL-map
onto $X$. When $X$ is a complete Riemannian manifold $M^{n}$ of dimension $%
n\geq 2$, the Menger sponge, the Sierpin'ski carpet or gasket, $\overline{X}$
is isometric to the so-called \textquotedblleft universal\textquotedblright\ 
$\mathbb{R}$-tree $A_{\mathfrak{c}}$, which has valency $\mathfrak{c=}%
2^{\aleph _{0}}$ at each point. In these cases, and when $X$ is the Hawaiian
earring $H$, the action of $\Gamma (X)$ on $\overline{X}$ gives examples in
addition to those of Dunwoody and Zastrow that negatively answer a question
of J. W. Morgan about group actions on $\mathbb{R}$-trees. Indeed, for one
length metric on $H$, we obtain precisely Zastrow's example.

MSC Classification: 57M07; 20F65, 28A80, 53C23, 54F15, 54F50
\end{abstract}

\keywords{length space, $\mathbb{R}$-tree, Hawaiian earring, fractal,
arcwise isometry, submetry, $\mathbb{R}$-free group}
\maketitle

\section{Introduction and Main Results}

A metric space such that every pair of its points is joined by a path of
length arbitrarily close to the distance between them is called an \textit{%
inner metric space} or \textit{length space}. A \textit{geodesic space} is a
metric space such that every pair of points is joined by a geodesic, i.e. a
path whose length is equal to the distance between them. Evidently every
geodesic space is a length space and it is a classical result that a
complete, locally compact length space is a geodesic space. A geodesic space
that contains no topological circle is called an $\mathbb{R}$-tree. A 
\textit{submetry} (resp. \textit{weak submetry}) $f:X\rightarrow Y$ between
metric spaces is a function which maps every closed (respectively, open)
ball in $X$ centered at any point $x\in X$ onto the closed (open) ball in $Y$
of the same radius at the point $f(x)$ (\cite{Bs}). Note that a submetry or
weak submetry is open, surjective and distance non-increasing, hence $1$%
-Lipschitz and uniformly continuous. A map is \textit{light} if every point
pre-image is totally disconnected (\cite{An}, \cite{W}). A function $%
f:X\rightarrow Y$ is a \textit{metric quotient} if $d_{Y}(x,y)$ is the
Hausdorff distance between $f^{-1}(x)$ and $f^{-1}(y)$; clearly a metric
quotient is a weak submetry. Recall that a group is \textit{locally free} if
each of its finitely generated subgroups is free.

\begin{theorem}
\label{ut1} Every length space (resp. complete length space) $(X,d)$ is the
metric quotient of a (resp. complete) $\mathbb{R}$-tree $(\overline{X},%
\overline{d})$ via the free isometric action of a locally free subgroup $%
\Gamma (X)$ of the isometry group $Isom(\overline{X})$. The quotient mapping 
$\overline{\phi }:\overline{X}\rightarrow X$ is a weak submetry (hence open)
and light map, and $\overline{\phi }$ is a submetry if $X$ is geodesic.
\end{theorem}

The $\mathbb{R}$-tree $\overline{X}$ is defined as the space of based
\textquotedblleft non-backtracking\textquotedblright\ rectifiable paths in $X
$, where the distance between two paths is the sum of their lengths from the
first bifurcation point to their endpoints. The group $\Gamma (X)\subset 
\overline{X}$ is the subset of loops with a natural group structure and the
quotient mapping $\overline{\phi }:\overline{X}\rightarrow X$ is the
end-point map. We will refer to $\overline{X}$ as the \textit{covering} $%
\mathbb{R}$\textit{-tree} of $X$. The term \textquotedblleft $\mathbb{R}$%
-tree\textquotedblright\ was coined by Morgan and Shalen (\cite{MS}) in 1984
to describe a type of space that was first defined by Tits (\cite{T}) in
1977. In the last three decades $\mathbb{R}$-trees have played a prominent
role in topology, geometry, and geometric group theory (see, for example, 
\cite{AM}, \cite{Be}, \cite{Ch}, \cite{MS}, \cite{Hu}). They are the most
simple of geodesic spaces, and yet Theorem \ref{ut1} shows that every length
space, no matter how complex, is an orbit space of an $\mathbb{R}$-tree.

Unless otherwise stated, \textquotedblleft dimension\textquotedblright\
refers to the covering dimension $\dim (X)$ of $X$. The small and large
inductive dimension of an a metric space $X$ satisfy the Katetov equality $%
Ind(X)=\dim (X)$ and the inequality $ind(X)\leq Ind(X)$ (see \cite{AP}). If $%
X$ is also separable, in particular compact, then $ind(X)=Ind(X)=\dim (X)$
(see \cite{HW}). The above theorem and the fact that a (non-trivial) $%
\mathbb{R}$-tree $X$ is simply connected with $ind(X)=1$ (\cite{AM}, \cite%
{Mo}, \cite{AB}) give us:

\begin{corollary}
\label{pc}Every non-trivial topological space admitting a compatible length
metric is the image via a light open mapping of a simply connected space $X$
with $ind(X)=1$.
\end{corollary}

Corollary \ref{pc} is broadly applicable because in 1949 Bing and Moise (%
\cite{B}, \cite{Moi}) independently and positively answered a 1928 question
of Menger: whether (in modern terminology) every Peano continuum (continuous
image of $[0,1]$) admits a compatible geodesic metric (\cite{M}). The
Bing-Moise Theorem shows that Corollary \ref{pc} contributes to a 70
year-old program to construct dimension-raising open mappings, beginning
with an example of Kolmogorov in 1937 (\cite{K}) from a Peano curve ($1$%
-dimensional Peano continuum) to a $2$-dimensional space. Later examples
include \cite{Ke} and the spectacular theorem partly stated without proof by
Anderson (\cite{An}, \cite{An2}) and proved by Wilson in \cite{W}: every
Peano continuum is the image via a light open mapping of the Menger sponge $%
\mathbb{M}$. Recall that $\mathbb{M}$ is called the \textquotedblleft
universal curve\textquotedblright\ because every Peano curve may be embedded
in it; the Anderson-Wilson Theorem provides a second sense in which $\mathbb{%
M}$ is \textquotedblleft universal\textquotedblright . We will be interested
in two more Peano curves: the Sierpin'ski carpet $S_{c}$ and gasket $S_{g}$.
As is well-known, each of these three \textit{fractal curves} (\cite{F})
admits a geodesic metric $d$ bi-Lipshitz equivalent to the metric induced by
the ambient Euclidean space of which it is a subspace: $d(x,y)$ is simply
the shortest Euclidean length of a path joining $x$ and $y$ in the space .

Next, recall that for a point $t$ in a $\mathbb{R}$-tree $T$, the valency at 
$t$ is the cardinality of the set of connected components of $T\backslash
\{t\}$, and $T$ is said to have valency at most $\mu $ if the valency of
every point in $T$ is at most $\mu $. A nontrivial complete metrically
homogeneous $\mathbb{R}$-tree can be characterized as a complete $\mathbb{R}$%
-tree $A_{\mu }$ with valency $\mu $ at each point for a cardinal number $%
\mu \geq 2$. It is unique up to isometry, and $\mu $-universal in the sense
that every $\mathbb{R}$-tree of valency at most $\mu $ isometrically embeds
in $A_{\mu }$. The existence of $A_{\mu }$ and the results just mentioned
were proved in \cite{MNO}. Another construction of $A_{\mu }$ was given in 
\cite{DP}, where it was shown that $A_{\mathfrak{c}}$ ($\mathfrak{c}%
=2^{\aleph _{0}}$, the cardinality of the continuum) can be isometrically
embedded at infinity in a complete simply connected Riemannian manifold of
constant negative curvature.

\begin{theorem}
\label{universal} If $X$ is a separable length space, then $\overline{X}$ is
a sub-tree of $A_{\mathfrak{c}}$. If in addition $X$ is complete and
contains a bi-Lipschitz copy of $S_{g}$ or $S_c$ at every point, e.g. if $X$
is $S_{c} $, $S_{g}$, $\mathbb{M}$, or a complete Riemannian manifold $M^n$
of dimension $n\geq 2$, then $\overline{X}$ is isometric to $A_{\mathfrak{c}%
} $.
\end{theorem}

Put another way, every separable length space may be obtained by starting
with a subtree of $A_{\mathfrak{c}}$ and taking a quotient of that subtree
via a free isometric action. Another consequence of this theorem is an
explicit construction of $A_{\mathfrak{c}}$ starting with any of the above
spaces (see the proof of Theorem \ref{ut1}). Notice that by using different
Banach spaces $X$ with their natural geodesic metric, we can in a similar
way realize $A_{\mu }$ as $\overline{X}$ for arbitrary $\mu \geq \mathfrak{c}
$. Our results, combined with the Anderson-Wilson Theorem, show that $A_{%
\mathfrak{c}}$ is \textquotedblleft universal\textquotedblright\ in a way
analogous to the second way in which $\mathbb{M}$ may be regarded as
\textquotedblleft universal\textquotedblright :

\begin{corollary}
Every non-trivial Peano continuum is the image of $A_{\mathfrak{c}}$ via a
light open mapping.
\end{corollary}

The function $\overline{\phi }$ from Theorem \ref{ut1} is generally not a
locally isometric covering map in the traditional sense, but shares
important properties with any such map $f:X\rightarrow Y$: (I) $f$ preserves
the length of rectifiable paths in the sense that $L(c)=L(f\circ c)$ for
every path $c$ in $X$ with finite length $L(c)$. (II) If $c$ is any
rectifiable path in $Y$ starting at a point $p$ and $f(q)=p$ then there is a
unique path $c_{L}$ starting at $q$ such that $f\circ c_{L}=c$, and moreover 
$c_{L}$ is rectifiable. A function $f$ between length spaces will be called 
\textit{unique rectifiable lifting (URL)} if it has these two properties.
Note that a map between length spaces with condition (I) is known as an 
\textit{arcwise isometry} (\cite{Gr}); such maps are a distant
generalization of isometric immersions from differential geometry. Any
URL-map is an arcwise isometric weak submetry (Proposition \ref{ws}). For
Riemannian manifolds, the notion of weak submetry is the same as that of
Riemannian submersion (\cite{BG}), which is in some sense dual to isometric
immersion. Generally a URL-map may not be locally injective at any point, as
Theorems \ref{ut1}, \ref{universal}, and \ref{m1} show.

\begin{theorem}
\label{m1}Under the assumptions and with the notation of Theorem \ref{ut1}:

\begin{enumerate}
\item The map $\overline{\phi }:\overline{X}\rightarrow X$ is a URL-map.

\item If $Z$ is a length space and $f:Z\rightarrow X$ is a URL-map then
there is a unique (up to basepoint choice) URL-map $\overline{f}:\overline{X}%
\rightarrow Z$ such that $\overline{\phi }=f\circ \overline{f}$.

\item There exists a unique (up to isometry) length space $(X_{1},d_{1})$
with a map $\phi _{1}:X_{1}\rightarrow X$, having the previous two
properties.

\item If there is an $\mathbb{R}$-tree $X_{1}$ with a URL-mapping $\phi
_{1}:X_{1}\rightarrow X$, then there is an isometry $\overline{\phi }%
_{1}:X_{1}\rightarrow \overline{X}$ such that $\phi _{1}=\overline{\phi }%
\circ \overline{\phi }_{1}$.
\end{enumerate}
\end{theorem}

In the language of category, this theorem means that $\overline{\phi }$ is
the initial object in the category of URL-mappings over $X$, i.e. $\overline{%
X}$ is \textquotedblleft universal\textquotedblright\ in this category. This
result, combined with Theorem \ref{universal}, shows that $A_{\mathfrak{c}}$
is \textquotedblleft universal\textquotedblright\ in yet a third sense.

One can easily deduce from Theorems \ref{ut1} and \ref{m1} the following
corollary.

\begin{corollary}
\label{samet} Let $f:(X_{1},\ast )\rightarrow (X_{2},\ast )$ be a basepoint
preserving URL-map of length spaces. Then there is a commutative diagram 
\begin{equation*}
\begin{array}{lll}
(\overline{X_{1}},\ast ) & \overset{\overline{f}}{\longrightarrow } & (%
\overline{X_{2}},\ast ) \\ 
\downarrow ^{\overline{\phi _{1}}} &  & \downarrow ^{\overline{\phi _{2}}}
\\ 
(X_{1},\ast ) & \overset{f}{\longrightarrow } & (X_{2},\ast )%
\end{array}%
\end{equation*}%
of URL-maps preserving basepoints, with unique $\overline{f}$, where $%
\overline{\phi _{1}}$ and $\overline{\phi _{2}}$ are the $\mathbb{R}$-tree
covering maps for $X_{1}$ and $X_{2}$ respectively, and $\overline{f}$ is
isometry. The identification $(\overline{X_{1}},\ast )$ with $(\overline{%
X_{2}},\ast )$ by the isometry $\overline{f}$ induces the homomorphic
inclusion $\Gamma (X_{1})\subset \Gamma (X_{2})$ of the corresponding
isometry groups.
\end{corollary}

Notice that the proof of Theorem \ref{ut1} implies that $\Gamma (X)$ is
naturally identified with the group $\Gamma $ from Proposition \ref{group}.
There naturally arises:

\begin{question}
For a given length space $(X,\ast)$, which subgroups of $\Gamma(X,\ast)$
correspond to URL-maps onto $(X,\ast)$?
\end{question}

We plan to occupy ourselves with this problem in the future. It may be
useful to consider the well-known bijective correspondence between
geodesically complete $\mathbb{R}$-trees with basepoints (\textquotedblleft
rooted $\mathbb{R}$-trees\textquotedblright ) and ultrametric spaces of
diameter 1 with nonempty spheres of radius 1 that come from considering the
end space of the $\mathbb{R}$-tree (\cite{Hu}). The ultrametric space
corresponding to $A_{\mu }$ is complete, metrically homogeneous, and does
not depend up to isometry on the choice of a basepoint in $A_{\mu }$. Then
an approach of Bruce Hughes in \cite{Hu} associates to each isometry of the $%
\mathbb{R}$-tree a so-called local similarity equivalence of the
corresponding ultrametric space. The preprint \cite{KW} may also be useful
for this problem.

\begin{remark}
Generally, for a given length space $X,$ we can find a proper subtree $%
\tilde{X}\subset \overline{X}$ and a proper subgroup $\tilde{\Gamma}%
(X)\subset \Gamma (X)$ such that $X$ is the metric quotient of the $\mathbb{R%
}$-tree $\tilde{X}$ via the free isometric action of the group $\tilde{\Gamma%
}(X)$ on $\tilde{X}$, and the quotient mapping $\tilde{\phi}:\tilde{X}%
\rightarrow X$ is an arcwise isometry, a weak submetry, and $\tilde{\phi}$
is a submetry if $X$ is geodesic. Under these conditions, $\tilde{X}$ is not
necessarily complete, even if $X$ is geodesic and complete. This is shown in
Theorem \ref{utr} for any Riemannian manifold $M^{n}$ of dimension $n\geq 2.$
It follows from Corollary \ref{samet} that $\tilde{\phi}$ is not a URL-map.
\end{remark}

Previously we discussed three main actors: any length space $X$, the $%
\mathbb{R}$-tree $\overline{X}$, and the URL-map $\overline{\phi }:\overline{%
X}\rightarrow X$. Theorem \ref{ut1} implies that for any pointed length
space $X$, the group $\Gamma (X)$ acts freely by isometries on the covering $%
\mathbb{R}$-tree $\overline{X}$. So our paper is closely connected with the
following general question of J. W. Morgan from \cite{Mo1}:

\begin{question}
\label{Mor} Which (finitely presented) groups act freely (by isometries) on $%
\mathbb{R}$-trees?
\end{question}

This question inspired us to study more closely the structure of the fourth
actor, the group $\Gamma (X)$. The answer to Question \ref{Mor} is known for
finitely generated groups (\cite{MS1}, \cite{GLP}, \cite{Ch}). However,
there are examples by Dunwoody \cite{Dun} and Zastrow \cite{Za} of
infinitely generated groups that are not free products of fundamental groups
of closed surfaces and abelian groups, but which act freely on an $\mathbb{R}
$-tree. Zastrow's group $G$ contains one of the two Dunwoody groups as a
subgroup. The other group is a Kurosh group. We prove the following theorems.

\begin{theorem}
\label{plenitude1} Let $X$ be $S_{c}$, $\mathbb{M}$, a complete Riemannian
manifold $M^{n}$ of dimension $n\geq 2$, or the Hawaiian earring $H$ with
any compatible length metric $d$. Then $\Gamma (X)$ is an infinitely
generated, locally free group that is not free and not a free product of
surface groups and abelian groups, but acts freely on the $\mathbb{R}$-tree $%
\overline{X}$. Moreover, the $\mathbb{R}$-tree $\overline{X}$ is a minimal
invariant subtree with respect to this action.
\end{theorem}

\begin{theorem}
\label{plenitude2} For any two length metrics $d_{1},d_{2}$ on $H$
(compatible with the usual topology), there is an injective homomorphism of $%
\Gamma (H,d_{1})$ into $\Gamma (H,d_{2})$. For a particular choice of $%
d=d_{Z}$ on $H$, $\Gamma (H,d_{Z})$ and its action on $\overline{(H,d_{Z})}$
coincide with Zastrow's group $G$ and its free action by isometries on
Zastrow's $\mathbb{R}$-tree.
\end{theorem}

An important role in the proofs is played by classical results about \textit{%
normal paths} from \cite{CF} and \textit{dendrites} from \cite{Ku}, and a
more recent characterization of $\mathbb{R}$-trees as Gromov \textit{%
0-hyperbolic} geodesic spaces (\cite{Gr}). There is an interesting
connection between these topics of different eras: It is not hard to show,
using the Bing-Moise Theorem, that a topological space is a dendrite if and
only if it is metrizable as a compact $\mathbb{R}$-tree.

This paper is connected with, and was inspired by, our previous paper \cite%
{BPU} and a 20 year old announcement of the first author, cited in \cite{DP}%
. The results from \cite{BPU} may be used to give an alternative proof of
Proposition \ref{rep}. In an upcoming paper we will show that \textit{all}
URL-maps, including the traditional universal cover of a length space, are
obtained via quotients of the covering $\mathbb{R}$-tree, and we will use
constructions of \cite{BPU} to obtain additional examples of URL-maps that
are not local isometries.

\textbf{Acknowledgements.} We thank Professors Karsten Grove, Ian Chiswell,
and Andreas Zastrow for very useful discussions. The first author is much
obliged to the Department of Mathematics of the University of Notre Dame in
Indiana for the hospitality he received as a visiting professor while part
of this paper was prepared. He was partially supported by RFBR (grant
08-01-00067-a) and the State Maintenance Program for the Leading Scientific
Schools of the Russian Federation (grant NSH-5682.2008.1).

\section{The covering $\mathbb{R}$-tree}

\begin{definition}
\label{wn} Let $c:[a,b]\rightarrow X$ be a path in a metric space $X$; $c$
is called normal if there is no nontrivial subsegment $J=[u,v]\subset
\lbrack a,b]$ such that $c(u)=c(v)$ and $c|_{J}$ is path homotopic to a
constant. Here \textquotedblleft path homotopic\textquotedblright\ means
fixed-endpoint homotopic. We define $c$ to be weakly normal if $c$ is normal
in its image $c([a,b])$.
\end{definition}

\begin{remark}
An immediate consequence of the above defintion is that every normal path is
weakly normal.
\end{remark}

\begin{proposition}
\label{ano}Consider the following three statements for a path $%
c:[a,b]\rightarrow X$:

\begin{enumerate}
\item $c$ is normal

\item $c$ is weakly normal

\item there is no nontrivial subsegment $J=[u,v]\subset \lbrack a,b]$ such
that $c(u)=c(v)$ and $c|_{J}$ is path-homotopic in $c(J)$ to a constant.
\end{enumerate}

All three are equivalent if $X$ is a separable, $1$-dimensional metric
space. If $X$ is an arbitrary metric space and $c$ is rectifiable then the
second and third are equivalent.
\end{proposition}

\begin{proof}
Suppose that $c$ is weakly normal and $X$ is separable and one-dimensional.
Then for every nontrivial subsegment $J=[u,v]\subset \lbrack a,b]$ such that 
$c(u)=c(v)=y$, $c|_{J}$ is not path homotopic in $Y:=c([a,b])$ to a constant
map, i.e. it represents a nontrivial element of $\pi _{1}(Y,y)$. By
Corollary 2.1 in \cite{CF}, the inclusion map $i:Y\rightarrow X$ induces an
injective homomorphism $i_{\ast }:\pi _{1}(Y,y)\rightarrow \pi _{1}(X,y)$.
So, the path $c|_{J}$ is not path-homotopic in $X$ to constant path. This
implies that the path $c$ is normal. The last statement follows rom the
previous statement and the well-known fact that the image $Z$ of a
nontrivial rectifiable path $c$ is one-dimensional. This fact follows from
inequalities $\dim (Z)\leq \dim _{H}(Z)\leq 1$, where $\dim _{H}$ is the
so-called Hausdorff dimension \cite{F}, \cite{Gr}, and non-triviality of $c$.
\end{proof}

\begin{definition}
\label{nonst} Two paths $c_{1},c_{2}:I=[a,b]\rightarrow X$ are called Fr\'{e}%
chet equivalent if there exist order-preserving monotone (continuous) maps $%
m_{1},m_{2}$ of $I$ onto itself such that $c_{1}m_{1}=c_{2}m_{2}$.
\end{definition}

In Lemma 3.1 and Theorem 3.1 the authors of \cite{CF} proved the following
results for a $1$\textit{-dimensional separable metric space} $X$.

\begin{lemma}
\label{normal} Each path $f:I\rightarrow X$ is path homotopic to a normal
path.
\end{lemma}

\begin{theorem}
\label{fresh} Two normal paths in $X$ are path homotopic if and only if they
are Fr\'{e}chet equivalent.
\end{theorem}

The definition of normal loop was given in \cite{CF} along with the
not-quite-standard Definition \ref{nonst}. The statements were proved for
loops, but the same arguments work for paths. Evidently, the
\textquotedblleft if\textquotedblright\ part of Theorem \ref{fresh} is valid
for any space $X$. Curtis and Fort proved Lemma \ref{normal} in \cite{CF} in
the following way. Let $S$ be the collection of all subsets $G$ of $I$, open
in $\mathbb{R}$, such that: If $(u,v)$ is a component of $G$, then $%
f(u)=f(v) $ and $f|_{[u,v]}$ is path homotopic to a constant. The collection 
$S$ is partially ordered by inclusion. It is proved that $S$ contains a
maximal element $G^{\ast }$. Define $g$ to be the map that agrees with $f$
on $I-G^{\ast }$ and is constant on each component of $G^{\ast }$. Then $g$
is path homotopic to $f$ and $g$ is normal.

\begin{proposition}
\label{rep} Any rectifiable path $c$ in a metric space is path homotopic in
its image to a weakly normal path $c_{n}$. Moreover, $L(c_{n})\leq L(c)$ and
the parameterization of $c_{n}$ by arclength is uniquely determined by $c$.
\end{proposition}

\begin{proof}
The first statement follows from Proposition \ref{ano} and Lemma \ref{normal}%
. The third statement is a corollary of the above hint for the proof of
Lemma \ref{normal}, Theorem \ref{fresh}, and the evident statement that two
rectifiable Fr\'{e}chet equivalent paths have equal parameterizations by
arclength.
\end{proof}

By a $\rho $\textit{-path} in a metric space $X$ we mean a weakly normal,
rectifiable, arclength parameterized path $c:[0,L]\rightarrow X$. Note that
the concatenation $c\ast d$ of a $\rho $-path $c$ followed by a $\rho $-path 
$d$ may not be a $\rho $-path. To resolve this problem we define the
\textquotedblleft cancelled concatenation\textquotedblright\ $c\star d$ to
be the unique $\rho $-path which is the arclength parameterization of the
weakly normal path in the path homotopy class of the concatenation $c\ast d$%
, in the image of $c\ast d$ (Proposition \ref{rep}). From the uniqueness and
the last statement in Proposition \ref{ano}, one can easily see more
concretely that $c\star d$ is obtained from $c\ast d$ by removing the
maximal final segment of $c$ that coincides with an initial segment of $d$
with reversed orientation, and removing that initial segment of $d$ as well.

\begin{proposition}
\label{group} The associative law $(a\star b)\star c=a\star (b\star c)$ is
satisfied. Moreover, cancelled concatenation on the set $\Gamma $ of all $%
\rho $-loops at a fixed basepoint $\ast $ of any metric space $X$ is a group
operation, where the constant loop is the identity and the inverse of $%
c:[0,L]\rightarrow X$ is the $\rho $-loop $c^{-1}(t):=c(L-t)$.
\end{proposition}

\begin{proof}
All these statements follow from the uniqueness of the $\rho$-path $c\star d$
for any two $\rho $-paths $c$ and $d$ in the case when $c\ast d$ makes sense
(see the discussion right before this proposition).
\end{proof}

We shall exclude in the future the trivial case when $X$ contains only one
point (this is traditional in discussing $\mathbb{R}$-trees). The following
are equivalent for a geodesic space $X$ (see \cite{Mo}, \cite{BH}, \cite{AB}%
): (1) $X$ is an $\mathbb{R}$-tree. (2) $X$ is $0$-hyperbolic in Gromov's
sense. (3) $X$ is $CAT(K)$-space for all $K\leq 0$. (4) $X$ is simply
connected and $ind(X)=1$ (\cite{AB}). See \cite{BH} for the definition of $%
CAT(K)$-space. We will not use this notion in the present paper except to
observe the corollary that every geodesic space is the metric quotient of a $%
CAT(K)$-space.

\begin{proof}[Proof of Theorem \protect\ref{ut1}]
Choose a base point $\ast \in X$ and define the set $\overline{X}$ to be the
set of all $\rho $-paths $c:[0,L]\rightarrow X$ starting at $\ast $. For $%
c_{1},c_{2}\in \overline{X}$, let $c_{1}\wedge c_{2}:[0,b]\rightarrow X$ be
the restriction of $c_{1}$ (and $c_{2}$) to the largest interval $[0,b]$ on
which $c_{1}$ and $c_{2}$ coincide, and define 
\begin{equation}
\overline{d}(c_{1},c_{2}):=L(c_{1})+L(c_{2})-2L(c_{1}\wedge
c_{2})=L(c_{1}^{-1}\star c_{2})\text{.}  \label{d}
\end{equation}

To see that $\overline{X}$ is an $\mathbb{R}$-tree, we will use the
characterization (2) above. We will also denote by $\ast $ the element of $%
\overline{X}$ that is simply the constant path at $\ast \in X$. Let $%
c_{1},c_{2}\in \overline{X}$, defined on $[0,L_{1}]$, $[0,L_{2}]$,
respectively. Let 
\begin{equation*}
s_{0}:=\max \{s:c_{1}(t)=c_{2}(t)\text{ for all }t\in \lbrack 0,s]\}
\end{equation*}%
and define $C(s)$ for $s\in \lbrack 0,L_{1}+L_{2}-2s_{0}]$ as follows. For $%
s\in \lbrack 0,L_{1}-s_{0}]$ let $C(s)$ be the restriction of $c_{1}$ to $%
[0,L_{1}-s]$. For $s\in \lbrack L_{1}-s_{0},L_{1}+L_{2}-2s_{0}]$ let $C(s)$
be the restriction of $c_{2}$ to $[0,s-L_{1}+2s_{0}]$. Certainly $C(s)$ is a
geodesic in $\overline{X}$ joining $c_{1}$ and $c_{2}$. This implies that $%
\overline{X}$ is a geodesic space.

We see from Formula (\ref{d}) that the so-called Gromov product 
\begin{equation*}
(c_{1},c_{2})_{\ast }:=\frac{1}{2}[\overline{d}(\ast ,c_{1})+\overline{d}%
(\ast ,c_{2})-\overline{d}(c_{1},c_{2})]
\end{equation*}%
with respect to the point $\ast $ (see, for example, \cite{BH}, page 410) is
equal to $L(c_{1}\wedge c_{2})$. Also we see immediately that $c_{1}\wedge
c_{2}$ contains as a subpath $(c_{1}\wedge c_{3})\wedge (c_{2}\wedge c_{3})$
for any $c_{3}$. Then it follows from these two statements that 
\begin{equation}
(c_{1},c_{2})_{\ast }\geq \min \{(c_{1},c_{3})_{\ast },(c_{3},c_{2})_{\ast
}\}  \label{c}
\end{equation}%
for any $c_{1},c_{2},c_{3}\in \overline{X}.$ This means that $\overline{X}$
is $0$-hyperbolic \textquotedblleft with respect to the point $\ast $%
\textquotedblright , whereas 0-hyperbolicity itself means that the Equation (%
\ref{c}) must be satisfied with respect to any point $c\in \overline{X}$.
But by Remark 1.21 (page 410 of \cite{BH}), $0$-hyperbolic at a single point
is sufficient for $\overline{X}$ to be $0$-hyperbolic, and hence an $\mathbb{%
R}$-tree.

By definition, $\overline{\phi }:(\overline{X},\overline{d})\rightarrow
(X,d) $ associates to a path $c\in \overline{X}$ its endpoint in $X$. Let $%
c_{1},c_{2}$ be any elements in $\overline{X}$. Then the path $%
c_{1}^{-1}\star c_{2}$ joins the points $x_{1}:=\overline{\phi }(c_{1})$ and 
$x_{2}:=\overline{\phi }(c_{s})$. Thus by definition 
\begin{equation}
d(x_{1},x_{2})\leq L(c_{1}^{-1}\star c_{2})=\overline{d}(c_{1},c_{2}).
\label{ine}
\end{equation}%
This means that $\overline{\phi }$ does not increase distances. Now let $%
x_{1},x_{2}\in X$, $x_{1}:=\overline{\phi }(c_{1})$ and $\varepsilon >0$ be
given. Then there is a rectifiable path $c$ in $X$ joinining the points $%
x_{1},x_{2}$ so that $L(c)\leq d(x_{1},x_{2})+\varepsilon $. Denote by $%
c_{n} $ the $\rho $-path in $X$ from Proposition \ref{rep}. Then $c_{n}$
joins the points $x_{1},x_{2}$, $L(c_{n})\leq L(c)$, and the $\rho $-path $%
c_{2}:=c_{1}\star c_{n}$ joins the points $\ast $ and $x_{2}$. Moreover, 
\begin{equation}
\overline{d}(c_{1},c_{2})=L(c_{1}^{-1}\star c_{2})=L(c_{1}^{-1}\star
c_{1}\star c_{n})=L(c_{n})\leq d(x_{1},x_{2})+\varepsilon \text{.}
\label{inee}
\end{equation}%
This together with Inequality (\ref{ine}) means that $\overline{\phi }$ is a
weak submetry, hence a metric quotient (\cite{Po}). If $(X,d)$ is a geodesic
space, we can take $c_{n}$ to be a geodesic and instead of the Inequality (%
\ref{inee}), we get 
\begin{equation*}
\overline{d}(c_{1},c_{2})=L(c_{1}^{-1}\star c_{2})=L(c_{1}^{-1}\star
c_{1}\star c_{n})=L(c_{n})=d(x_{1},x_{2})\text{.}
\end{equation*}%
This together with Inequality (\ref{ine}) means that $\overline{\phi }$ is a
submetry.

For an arbitrary interior point $w$ of a non-trivial geodesic segment $[y,z]$
in an $\mathbb{R}$-tree $T$, $y,z$ lie in different connected components of $%
T-\{w\}$. Then any connected subset $C\subset \overline{X}$ containing two
different points $y,z$ must include $[y,z]$. Now it follows from the
definition that in this case $\overline{\phi }([y,z])$ is the image of a
non-trivial path in $X$. Then the inclusion $C\subset \overline{\phi }%
^{-1}(x)$ is impossible for any $x\in X$, which shows that $\overline{\phi }$
is a light map.

Now let $\Gamma $ be the group from Proposition \ref{group}. It follows from
Formula (\ref{d}) and Proposition \ref{group} that $\Gamma $ acts freely via
isometries on $(\overline{X},\overline{d})$ if we define $l(c)=l\star c$ for
any $\rho $-loop $l\in \Gamma $ and $c\in \overline{X}$. Obviously, $%
\overline{\phi }(l(c))=\overline{\phi }(c)$. Also, if $c_{1},c_{2}\in 
\overline{X}$ and $c_{1},c_{2}\in \overline{\phi }^{-1}(x)$, $x\in X$, then $%
c_{2}=l(c_{1})$, where $l=c_{2}\star c_{1}^{-1}$. This means that every
pre-image $\overline{\phi }^{-1}(x)$, $x\in X$, is an orbit via the action
of $\Gamma $. The orbits of this action are precisely the sets $\phi
^{-1}(x) $ for $x\in X$, which verifies that $X$ is the metric quotient with
respect to the action of $\Gamma $.

Assume that $X$ is complete. Suppose that $c_{k}:[0,L_{k}]\rightarrow X$ is
a Cauchy sequence in $X$. By definition of the metric, $\{L_{k}\}$ converges
to a real number $L$ and so is bounded above by some finite number $M$. By
extending all paths to be constant at their endpoints we may assume that all
paths are defined on $[0,M]$ (these extensions generally are not in $%
\overline{X}$). Now all these paths are 1-Lipshitz maps. That is, the
sequence of these paths is uniformly Cauchy and since $X$ is complete, it
converges uniformly to some path $c:[0,M]\rightarrow X$. It follows from the
uniform convergence that $c_{[0,L]}\in \overline{X}$ and $c_{k}\rightarrow
c_{[0,L]}$ in $(\overline{X},\overline{d})$. This proves the completeness of 
$\left( \overline{X},\overline{d}\right) $.

Finally we check that $\Gamma (X)$ is locally free. Let $l_{1},...,l_{n}$ be 
$\rho $-loops in $X$ starting at $\ast $ and $Y$ be the union of their
images, which is separable and $1$-dimensional. According to \cite{CF1}, $%
\pi _{1}(Y)$ is locally free. Evidently the subgroup $\Gamma
(l_{1},...,l_{n})$ of $\Gamma (X)$ generated by $l_{1},...,l_{n}$ is
naturally identified with a subgroup of $\Gamma (Y)$. Moreover, it follows
from the definition of $\Gamma (X)$, Lemma \ref{normal}, and Theorem \ref%
{fresh} that $\Gamma (Y)$ is naturally isomorphic to a subgroup of the
locally free group $\pi _{1}(Y)$. By the Nielsen-Schreier Theorem, $\Gamma
(l_{1},...,l_{n})$ is free.
\end{proof}

\begin{lemma}
\label{geod} Let $\gamma _{c}$ denote the unique geodesic in $\overline{X}$
parameterized by arclength and joining the points $\ast $ (constant path at
the point $\ast \in X$) and $c$. Then

\begin{enumerate}
\item $\overline{\phi }\circ \gamma _{c}=c$,

\item $\gamma _{c_{1}}^{-1}\star \gamma _{c_{2}}$ is the unique geodesic in $%
\overline{X}$ parameterized by arclength and joining the points $c_{1}$ and $%
c_{2}$,

\item $\overline{\phi }\circ (\gamma _{c_{1}}^{-1}\star \gamma
_{c_{2}})=c_{1}^{-1}\star c_{2}$, and

\item $L(\gamma _{c_{1}}^{-1}\star \gamma _{c_{2}})=L(c_{1}^{-1}\star c_{2})$%
.
\end{enumerate}
\end{lemma}

\begin{proof}
These statements follow from definition of $\overline{d}$, the fact that the
path $C$ considered in the proof of Theorem \ref{ut1} is the unique
arclength-parameterized geodesic in $(\overline{X},\overline{d})$ joining
the points $c_{1},c_{2}\in \overline{X}$, and the equations $C=\gamma
_{c_{1}}^{-1}\star \gamma _{c_{2}}$, $\overline{\phi }\circ
C=c_{1}^{-1}\star c_{2}$.
\end{proof}

In fact, the construction of $\overline{X}$ and $\Gamma (X)$ in the proof of
Theorem \ref{ut1} depends on the choice of the base point $\ast \in X$, so,
strictly speaking, we must write $\overline{(X,\ast )}$ and $\Gamma (X,\ast )
$ instead of $\overline{X}$ and $\Gamma (X)$. Proposition \ref{equivo} below
shows that this dependence is not so essential. Recall the following
definition.

\begin{definition}
\label{equiv}An action of a group $\Gamma _{1}$ on a metric space $X_{1}$
via isometries is said to be equivalent to an action of a group $\Gamma _{2}$
on a metric space $X_{2}$ via isometries if there are an isometry $%
f:X_{1}\rightarrow X_{2}$ and an isomorphism $\phi :\Gamma _{1}\rightarrow
\Gamma _{2}$ such that $f(g(x))=\phi (g)(f(x))$ for any point $x\in X$ and
any element $g\in \Gamma _{1}$.
\end{definition}

\begin{proposition}
\label{equivo}Let $X$ be a length space and $\ast ,\star \in X$ be two of
its points. Then the action of the group $\Gamma (X,\ast )$ on the $\mathbb{R%
}$-tree $\overline{(X,\ast )}$ is equivalent to the action of the group $%
\Gamma (X,\star )$ on the $\mathbb{R}$-tree $\overline{(X,\star )}$.
\end{proposition}

\begin{proof}
By Proposition \ref{rep}, there is a $\rho $-path $k$ in $X$ which starts at 
$\star $ and ends at $\ast $. Define maps $f:\overline{(X,\ast )}\rightarrow 
\overline{(X,\star )}$ and $\phi :\Gamma (X,\ast )\rightarrow \Gamma
(X,\star )$ by the formulas $f(c):=k\star c$ and $\phi (\gamma )=k\star
\gamma \star k^{-1}$. By Proposition \ref{group} and Formula (\ref{d}), 
\begin{equation*}
\overline{d}(f(c_{1}),f(c_{2}))=\overline{d}(k\star c_{1},k\star
c_{2})=L((k\star c_{1})^{-1}\star (k\star c_{2}))
\end{equation*}%
\begin{equation*}
=L(c_{1}^{-1}\star k^{-1}\star k\star c_{2})=L(c_{1}^{-1}\star c_{2})=%
\overline{d}(c_{1},c_{2})
\end{equation*}%
for any two $\rho $-paths $c_{1},c_{2}\in \overline{(X,\ast )}$; for any
element $c^{\prime }\in \overline{(X,\star )}$, $k^{-1}\star c^{\prime
}:=c\in \overline{(X,\ast )}$ and $f(c)=k\star k^{-1}\star c^{\prime
}=c^{\prime }$. Thus $f$ is an isometry. Now by Proposition \ref{group}, 
\begin{equation*}
\phi (\gamma _{1}\star \gamma _{2})=k\star (\gamma _{1}\star \gamma
_{2})\star k^{-1}=(k\star \gamma _{1}\star k^{-1})\star (k\star \gamma
_{2}\star k^{-1})=\phi (\gamma _{1})\star \phi (\gamma _{2})
\end{equation*}%
for any $\gamma _{1},\gamma _{2}\in \Gamma (X,\ast )$, and for any element $%
\gamma ^{\prime }\in \Gamma (X,\star )$, 
\begin{equation*}
\gamma ^{\prime -1}\star \gamma ^{\prime -1}=\phi (k^{-1}\star \gamma
^{\prime }\star k)\text{, where }k^{-1}\star \gamma ^{\prime }\star k\in
\Gamma (X,\ast )\text{.}
\end{equation*}%
So $\phi $ is an isomorphism. Finally, 
\begin{equation*}
f(\gamma (c))=k\star (\gamma \star c)=(k\star \gamma \star k^{-1})\star
(k\star c)=\phi (\gamma )(f(c))
\end{equation*}%
for any elements $c\in \overline{(X,\ast )}$ and $\gamma \in \Gamma (X,\ast )
$.
\end{proof}

\section{Continua, Fractals, and Manifolds}

For the next proposition, let $X$ be a length space, $p\in X$. Define $\rho $%
-paths $\alpha :[0,a]\rightarrow X$ and $\beta :[0,b]\rightarrow X$ starting
at $p$ to be equivalent if $\alpha ,\beta $ coincide on $[0,\varepsilon )$
for some $\varepsilon >0$. We denote the cardinality of the set of the
resulting equivalence classes by $\kappa _{p}$.

\begin{proposition}
\label{val}The valency of $\overline{X}$ at any point $\overline{p}\in 
\overline{\phi }^{-1}(p)$ is equal to $\kappa _{p}$. If $X$ is separable
then $\kappa _{p}\leq \mathfrak{c}=2^{\aleph _{0}}$.
\end{proposition}

\begin{proof}
The construction of $\overline{X}$ immediately implies the first statement.
If $X$ is separable, then $X$ itself has cardinality $\mathfrak{c}$ (unless
it is a point). Since every path is determined by its value at rational
numbers in its domain, the cardinality of $\kappa _{p}$ is at most $%
(2^{\aleph _{0}})^{\aleph _{0}}=2^{\aleph _{0}\cdot \aleph _{0}}=2^{\aleph
_{0}}=\mathfrak{c}$.
\end{proof}

\begin{example}
\label{Zas}Consider the space $B$ consisting of countably many circles $%
\left\{ C_{i}\right\} _{i\in \mathbb{N}}$ each of length $s_{i}>0$, all
attached at a common basepoint $\ast $, and given the induced geodesic
metric. We are interested in two cases: (1) $s_{i}$ is a constant $s$; (2) $%
s_{i}$ is a strictly decreasing sequence converging to zero. In the first
case $B$ is not compact. The valency of points in $\overline{B}$ is either $%
\aleph _{0}$ or $2$ depending on whether they are in $\overline{\phi }%
^{-1}(p)$ or not. In the second case $B$ is homeomorphic to the Hawaiian
earring $H$, considered below.
\end{example}

An early result concerning $H$ is Theorem 2.1 from \cite{CF} which states:
If a one-dimensional separable locally connected continuum is not locally
simply connected, then it contains a subspace which has the homotopy type of 
$H$.

\begin{proposition}
\label{haw}For the Hawaiian earring $H$, supplied with any length metric of
the second kind from Example \ref{Zas}, we have $\kappa _{p}=\mathfrak{c}$
for the point $p=\ast $ and $\kappa _{p}=2$ for any point $p\neq \ast $.
\end{proposition}

\begin{proof}
The last statement is evident. It is easy to find $\mathfrak{c}$ unit weakly
normal loops starting (and ending) at $\ast $ such that no two coincide on
any interval $[0,\varepsilon )$. In fact, let us take first an increasing
integer sequence $\{i(n),n\in \mathbb{N}\}$ so that $s_{i(n)}<\frac{1}{2^{n}}%
.$ One can define a path that wraps one of two ways around $C_{i(1)}$, then
one of two ways around $C_{i(2)}$, and so on. Then reverse the
parametrization, so that $C_{i(1)}$ is wrapped around last. It is clear that
any such path is weakly normal and is encoded as a sequence $\{x_{n}\}_{n\in 
\mathbb{N}}$ with values in the set $\{-1,1\}$. Two such arclength
parameterized paths are equivalent if and only if they wrap the same way
around $C_{i(n)}$ for all sufficiently large $n$, or in other words, define
equivalent sequences $\{x_{n}\}$ and $\{y_{n}\}$--we mean here that $x$ is
equivalent to $y$ if there is a number $m\in \mathbb{N}$ such that $%
x_{n}=y_{n}$ for all $n\geq m$. The following lemma finishes the proof.
\end{proof}

\begin{lemma}
We have $\mathfrak{c}$ different equivalence classes of the above type of
sequences.
\end{lemma}

\begin{proof}
Let us take first the sequence $x_{n}\equiv 1$. Consider now an arbitrary
sequence $z=\{z_{n}\}$ of natural numbers and define subsequently two
sequences $s(z)_{n}:=\sum_{i=1}^{n}z_{i}$ and $\sigma
(z)_{n}:=\sum_{i=1}^{n}2^{s(z)_{i}}$. It is clear that if for another such
sequence $w=\{w_{n}\}$, $w_{m}\neq z_{m}$, then $\sigma (w)_{n}\neq \sigma
(z)_{n}$ for all $n\geq m$. Then for every sequence $z$ above, define
another sequence $x(z)_{n}$ with values in $\{-1,1\}$ by the equations $%
x(z)_{m}=-1$ if $m=\sigma (z)_{k}$ for some $k\in \mathbb{N}$, and $%
x(z)_{m}=1$ for all other $m\in \mathbb{N}$. Then it follows from the
statement above that $x(z)_{n}$ is not equivalent to $x$ and is not
equivalent to $x(w)$ if $w\neq z$. So we get $\aleph _{0}^{\aleph _{0}}=%
\mathfrak{c}$ pairwise non-equivalent sequences of above sort.
\end{proof}

\begin{proof}[Proof of Theorem \protect\ref{universal}]
The first part of the theorem is an immediate consequence of Proposition \ref%
{val} and the theorem from \cite{MNO} that every $\mathbb{R}$-tree of
valency at most $\mathfrak{c}$ isometrically embeds in $A_{\mathfrak{c}}$.
Next, for any point $p\in S_{c}$ there is clearly a bi-Lipschitz embedding $%
h:H\rightarrow S_{c}$ such that $h(\ast )=p$. Therefore $\kappa _{p}\geq 
\mathfrak{c}$. The case of $S_{g}$ is more tricky. There are countably many
rectifiable loops $C_{i}$ in $S_{g}$ starting at any fixed point $p$ such
that every $C_{i}$ is a topologically embedded circle and $L(C_{i+1})<\frac{1%
}{3}L(C_{i})$ for all natural numbers $i$. Then one can prove with a little
more detail than in Proposition \ref{haw}, that $\kappa _{p}\geq \mathfrak{c}
$. A similar argument now finishes the proof of the theorem.
\end{proof}

\section{URL-maps}

The following statement easily follows from definitions.

\begin{proposition}
\label{cat} The collection of pointed length spaces with URLs as morphisms
is a category.
\end{proposition}

\begin{proposition}
\label{mc} A map $f:X\rightarrow Y$ between length spaces is a URL-map if
and only if $f$ is 1-Lipschitz and for some choice of basepoints, $f$ is
basepoint preserving, each arclength parameterized rectifiable path $p$ in $%
Y $ starting at the basepoint has unique lift $p_{L}$ starting at the
basepoint, and $L(p)=L(p_{L})$.
\end{proposition}

\begin{proof}
The necessity of these conditions easily follows from the definition of
URL-map. Let us prove sufficiency. Assume first that $c$ is an arclength
parameterized rectifiable path starting at $y\in Y$ and $x\in X$ satisfies $%
f(x)=y$. Let $k$ be a $\rho $ -path from the basepoint in $X$ to $x$. Then $%
d:=f\circ k$ is rectifiable and since each of its initial segments has, by
assumption, a unique lift of the same length, $d$ is also arclength
parameterized and has the same length as $k$. Then $d\ast c$ is rectifiable
and arclength parameterized, so has a unique lift $(d\ast c)_{L}$ to the
basepoint in $X$. We may write $(d\ast c)_{L}=k\ast k^{\prime }$ for some
path $k^{\prime }$. Then $k^{\prime }$ is the desired lift of $c$; $%
k^{\prime }$ must be unique since if it were not then $d\ast c$ would not
have a unique lift. Now suppose that $c:[0,a]\rightarrow Y$ is a rectifiable
path starting at $y$, and $f(x)=y$. Let $\mathcal{C}$ be the collection of
maximal (closed) intervals on which $c$ is constant. As is well-known, there
are a non-decreasing continuous function $h:[0,a]\rightarrow \lbrack 0,L(c)]$
and an arclength parameterized rectifiable path $c_{1}:[0,L(c)]\rightarrow Y$
such that $c=c_{1}\circ h$ and $h$ is strictly increasing everywhere except
on the intervals in $\mathcal{C}$. Let $d_{1}$ be the unique lift of $c_{1}$
at $x$. Define $c_{L}:[0,a]\rightarrow X$ by $c_{L}(t)=d_{1}\circ h(t)$.
Then $c_{L}$ has the same length as $d_{1}$, hence as $c_{1}$ and $c$. Since 
$f$ is 1-Lipshitz, it now follows from the lifting property proved above
that if $c$ is rectifiable in $X$ then $f\circ c$ has the same length as $c$%
. If $c$ is not rectifiable then $f\circ c$ cannot be rectifiable either,
for if it were, $f\circ c$ would have a rectifiable lift and a
non-rectifiable lift.
\end{proof}

\begin{proposition}
\label{ws} Every URL-map $f:X\rightarrow Y$ is a weak submetry. If $Y$ is a
geodesic space then $f$ is a submetry.
\end{proposition}

\begin{proof}
Since $Y$ is a length space, for any $\varepsilon >0$ and $x,y\in Y$ we may
join $x,y$ by a rectifiable path $c$ with the length less than $%
d(x,y)+\varepsilon $. By definition, $c$ has a lift $c_{L}$ of the same
length with endpoints $w,z$ such that $f(w)=x$ and $f(z)=y$. Since $X$ is a
length space, $d(w,z)\leq L(c_{L})=L(c)\leq d(x,y)+\varepsilon $ and the
proof of the first part is complete. If $Y$ happens to be a geodesic space
then we may take $c$ to be a geodesic with corresponding inequality $%
d(w,z)\leq L(c_{L})=L(c)=d(x,y)$. But from the first part we have that $f$
is 1-Lipschitz, so $d(x,y)\leq d(z,w)$ and $d(x,y)=d(z,w)$. This implies
that $f$ is a submetry.
\end{proof}

\begin{proof}[Proof of Theorem \protect\ref{m1}]
(1) By Theorem \ref{ut1}, $\overline{\phi }$ is a weak submetry, so is a
1-Lipschitz map. Now, by Proposition \ref{mc}, we need only to prove
condition (II) for points $p=\ast \in X$ and $q=\ast \in \overline{X}$ and
rectifiable arclength parameterized paths starting at $\ast \in X$. To do
so, let $c:[0,L]\rightarrow \overline{X}$ be any rectifiable path in $(X,d)$
starting at $\ast .$ By Proposition \ref{rep} we have, for every $0\leq
s\leq L$, a unique up to Fr\'{e}che equivalence weakly normal path $%
p_{s}:[0,s]\rightarrow c([0,s])$ such that $p_{s}$ is path homotopic to $%
c_{s}:=c\mid _{\lbrack 0,s]}$ in $c([0,s])$ and $L(p_{s})\leq L(c_{s})$.
Define a $\rho $-path in $(X,d)$ (and so an element in $\overline{X}$) $%
\overline{c}(s)$ to be the arclength parameterization of $p_{s}$.
Proposition \ref{rep} implies that $\overline{c}(s)$ is uniquely defined by $%
c$ and $s\in \lbrack 0,L]$. If $0\leq s_{1}<s_{2}\leq L$ then by definition
of $\overline{d}$, 
\begin{equation*}
\overline{d}(\overline{c}(s_{1}),\overline{c}(s_{2}))=L(\overline{c}%
(s_{1})^{-1}\star \overline{c}(s_{2}))\text{.}
\end{equation*}%
It is clear that the $\rho $-path $\overline{c}(s_{1})^{-1}\star \overline{c}%
(s_{2})$ is similarly defined by the path $c_{s_{1}}^{-1}\ast
c_{s_{2}}=c_{s_{1},s_{2}}:=c\mid _{\lbrack s_{1},s_{2}]}$. It follows from
Proposition \ref{rep} and the argument above that $\overline{d}(\overline{c}%
(s_{1}),\overline{c}(s_{2})\leq L(c_{s_{1},s_{2}})$. This implies that the
path $\overline{c}(s)$, $s\in \lbrack a,b]$ is continuous in $(\overline{X},%
\overline{d})$, $\overline{\phi }\circ \overline{c}=c$, and $L(\overline{c}%
)\leq L(c)$. Finally, since $\overline{\phi }$ is a 1-Lipschitz map we have $%
L(\overline{c})=L(c)$.

To finish the proof of the condition (II), we need to prove that if $%
\overline{c}^{\prime }:[0,L]\rightarrow \overline{X}$ is any path such that $%
\overline{c}^{\prime }(0)=\ast $ and $\overline{\phi }\circ \overline{c}%
^{\prime }=c$ then $\overline{c}^{\prime }=\overline{c}$. Since $\overline{X}
$ is an $\mathbb{R}$-tree by Theorem \ref{ut1}, for any $s\in \lbrack 0,L]$, 
$C_{s}=\overline{c}^{\prime }([0,s])$ is a Peano continuum that contains no
topological circle. By the Hahn-Mazurkiewicz-Sierpin'ski Theorem, $C_{s}$ is
locally connected, hence a dendrite (see Section 51, VI of \cite{Ku}). Then $%
C_{s}$ is contractible by Corollary 7 in Section 54, VII of \cite{Ku}. Hence
there is a path homotopy $h_{s}:[a,s]\times \lbrack 0,1]\rightarrow C_{s}$
such that $h(\cdot ,0)=\overline{c}^{\prime }|_{[0,s]}$, and $h(\cdot ,1):=%
\overline{c}_{s}^{\prime }$ is a topological embedding whose image is the
unique arc $a_{s}$ in $C_{s}$, joining $\overline{c}^{\prime }(0)=\ast $ and 
$\overline{c}^{\prime }(s)$ (this arc exists by Corollary 7, Section 51, VI
in \cite{Ku}). Since $\overline{X}$ is an $\mathbb{R}$-tree, the arc $a_{s}$
is a geodesic segment in $(\overline{X},\overline{d})$. By Lemma \ref{geod}, 
$\overline{\phi }\circ \overline{c}_{s}^{\prime }=\overline{\phi }\circ
\gamma _{\overline{c}^{\prime }(s)}=\overline{c}^{\prime }(s)$. It is clear
that $\overline{\phi }\circ h_{s}$ is a path homotopy in $c([0,s])$ from $%
c_{s}:=c|_{[0,s]}$ to the path $\overline{c}^{\prime }(s)$. So the last path
coincides with the $\rho $-path $\overline{c}(s)$ considered above, and we
have proved the required equality $\overline{c}^{\prime }=\overline{c}$.

(2). Given a URL-map $f:Z\rightarrow X$ with some choice of basepoints
define $\overline{f}(c)$ to be the endpoint of the unique lift of $c$
starting at the basepoint in $Z$. Obviously $f\circ \overline{f}=\overline{%
\phi }$. For $c,k\in \overline{X}$, the lift of $c\star k^{-1}$ is a path
joining $f(c)$ and $f(k)$ having the same length as $c\star k^{-1}=d(c,k)$,
and therefore $\overline{f}$ is $1$-Lipschitz. Now let $\gamma $ be a
rectifiable path starting at the basepoint in $Z$. Then $f\circ \gamma $ is
rectifiable in $X$, so has a lift $(f\circ \gamma )_{L}$ at the basepoint in 
$\overline{X}$. Now $\overline{f}\circ (f\circ \gamma )_{L}$ is a lift of $%
f\circ \gamma $ starting at the basepint in $Z$ and so must be equal to $%
\gamma $. That is, $(f\circ \gamma )_{L}$ is a lift of $\gamma $ starting at
the basepoint in $\overline{X}$ having the same length as $\gamma $. Suppose 
$k$ is any lift of $\gamma $ starting at the basepoint in $\overline{X}$.
Then $\gamma $ is also a lift of $f\circ \gamma $ to $\overline{X}$ and
therefore $k=(f\circ \gamma )_{L}$. We have checked the conditions of
Proposition \ref{mc} to show that $\overline{f}$ is a URL-map. Finally,
suppose we have a URL-map $h$ that preserves the basepoints with $f\circ h=%
\overline{\phi }$. For any $c\in \overline{X}$, $h\circ \gamma _{c}$ is a
rectifiable path from the basepoint to $h(c)$, which is also a lift of $%
\overline{\phi }\circ \gamma _{c}=c$ starting at the basepoint in $Z$. Since
this lift is unique, $h(c)=\overline{f}(c)$.

(3). The uniqueness of $\overline{X}$ follows from Proposition \ref{ws} and
the second part of the theorem.

(4). Let $(X_{1},d_{1})$ be an $\mathbb{R}$-tree and $\phi
_{1}:(X_{1},d_{1})\rightarrow (X,d)$ be a URL-map that preserves some
basepoints $\ast $. By Property (2) there is a unique basepoint preserving
URL-map $\overline{\phi _{1}}:\overline{X}\rightarrow X_{1}$ such that $%
\overline{\phi }=\phi _{1}\circ \overline{\phi _{1}}$. Since $X_{1}$
contains no topological circle, using the same arguments as in the proof of
condition (II) above, we get that any $\rho $-path $c$ in $X_{1}$ is
injective (otherwise we could prove that $c$ is not weakly normal), hence a
topological embedding and a geodesic in $X_{1}$. This together with the
construction of $\overline{X_{1}}$ implies that we can take $\overline{X_{1}}%
=X_{1}$ and $\overline{\phi }_{1}=id_{X_{1}}$. Again applying Property (2),
we find a URL-map $\overline{\overline{\phi _{1}}}:X_{1}\rightarrow 
\overline{X}$ such that $\overline{\phi _{1}}\circ \overline{\overline{\phi
_{1}}}=id_{X_{1}}$. Since all three of these maps are weak submetries by
Proposition \ref{ws}, they are all isometries.
\end{proof}

Using the same argument as in the proof for Part (4) of Theorem \ref{m1}, we
get

\begin{proposition}
If $f:X\rightarrow Y$ is URL-map between length spaces, and $Y$ is an $%
\mathbb{R}$-tree, then $f$ is an isometry.
\end{proposition}

\section{$\mathbb{R}$-free groups}

The primary reference for the discussion that follows is \cite{Ch}. A group
acting freely on an $\mathbb{R}$-tree is usually called $\mathbb{R}$\textit{%
-free}.

\begin{theorem}
(\cite{MS1})\label{MS1}The fundamental group of a closed surface is $\mathbb{%
R}$-free, except for the non-orientable surfaces of genus 1, 2, and 3 (the
connected sum of 1, 2, or 3 real projective planes).
\end{theorem}

The non-orientable surfaces of genus 1, 2, and 3 are called \textit{%
exceptional} and their fundamental groups \textit{exceptional surface groups}%
. The torus has fundamental group $\mathbb{Z}\oplus \mathbb{Z}$ and is
embeddable in $(\mathbb{R},+)$ (\cite{Ch}). Any subgroup of $(\mathbb{R},+)$
acts freely by isometries on $\mathbb{R}$, so is $\mathbb{R}$-free.

\begin{question}
(\cite{Mo})\label{Mo1}It follows easily that any free product of
non-exceptional surface groups and subgroups of $(\mathbb{R},+)$ is $\mathbb{%
R}$-free. The question is whether the converse statement is true.
\end{question}

The positive answer in the case of finitely generated groups is given by the
following \textit{Rips' Theorem}:

\begin{theorem}
(\cite{GLP}, \cite{BF}, \cite{Ch})\label{Rips} Any finitely generated $%
\mathbb{R}$-free group $G$ can be written as a free product $G=G_{1}\ast
\cdot \cdot \cdot \ast G_{n}$ for some integer $n\geq 1$, where each $G_{i}$
is either a finitely generated free abelian group or a non-exceptional
surface group.
\end{theorem}

As was pointed out in the Introduction, the answer to Question \ref{Mo1} is
negative in general. All spaces below are length spaces with basepoints and
maps are basepoint-preserving.

\begin{theorem}
\label{lynd}Let $L(g)$ denote the length of an element $g\in G=\Gamma
(X,\ast )$, or in other words, $L(g)=\overline{d}(g,\ast )$ if we consider $%
g $ as an element of $\overline{X}$. Then

\begin{enumerate}
\item $L(g)\geq 0$, and $L(g)=0$ if and only if $g=1=\ast $.

\item For all $g\in G$, $L(g)=L(g^{-1})$.

\item For all $g,h,k\in G,c(g,h)\geq \min \{c(g,k),c(h,k)\}$, where $c(g,h)$
is defined to be $\frac{1}{2}(L(g)+L(h)-L(g^{-1}h))$.
\end{enumerate}
\end{theorem}

\begin{proof}
The first two properties are evident. The third statement is simply the $0$%
-hyperbolicity property of $(\overline{X},\overline{d})$ for elements $g,h,k$%
, since we earlier referred to $c(g,h)$ as the Gromov product.
\end{proof}

This theorem implies that $L$ is a \textit{Lyndon length function} on $G$,
since by definition such a function satisfies precisely these properties
except that in condition (1) only the \textquotedblleft
if\textquotedblright\ part of the second statement is required. We will
refer to a Lyndon function satisfying the stronger condition (1) as \textit{%
definite}.

Now let $L$ be any definite Lyndon function on a group $G$. To obtain an $%
\mathbb{R}$-tree $T(G,L)$, first join any two elements $g,h\in G$ by an edge 
$[g,h]$ of length $L(g^{-1}h)$; then for any three elements $g,h,k\in G$
isometrically glue the initial segments of length $c(k^{-1}g,k^{-1}h)$ of
the edges $[k,g]$ and $[k,h]$ starting at $k$. By construction, any point $%
x\in T(G,L)$ lies in some edge $[g,h]$ for $g,h\in G$. Since any two
elements of $G$ are joined by a segment, any two points $x,y\in T(G,L)$ are
joined by some segment (of a finite length). Really this segment $[x,y]$ is
unique, and we can define $\rho (x,y)$ as the length of the segment $[x,y]$.
Now the action of $G$ on itself by the left multiplication, defined by
formula $l(g)(h)=gh$, has a well-defined extension to $T(G,L)$ by the
requirement that any segment $[h,k]$ maps isometrically onto the segment $%
[gh,gk]$. This is possible because $L((gh)^{-1}(gk))=L(h^{-1}k)$. Then $G$
acts freely on $T(G,L)$ by isometries. It is clear that $(T(G,L),\rho )$ and
the action of $G$ on $T(G,L)$ are uniquely defined by the function $L$ and
the described construction. The equality $L(g)=\rho (1,g)$ returns us to the
initial function $L$. Of course we have omitted the explanations of some
natural questions arising in the process of this construction, but these
details may be found in the literature on the subject.

Suppose now that $G$ acts freely by isometries on an $\mathbb{R}$-tree $%
(T,\rho )$. Choose any point $x\in T$ and define $L_{x}(g)=\rho (x,g(x))$.
Then $L_{x}$ is a Lyndon function on $G$ that in general depends on $x$. We
shall identify an element $g\in G$ with the point $g(x)\in T$. The $\mathbb{R%
}$-tree $T$ is isometric to $T(G,L_{x})$ and the action of $G$ on $T$ is
equivalent to the action of $G$ on $T(G,L_{x})$ if and only if $T$ is a 
\textit{minimal} $\mathbb{R}$-tree in $T$ containing all the elements of $G$.

\begin{proposition}
\label{Za}Zastrow's group $G=G_{Z}$ coincides with $\Gamma (B)$ for the
space $B$ from Example \ref{Zas} corresponding to the sequence $\{s_{i}=%
\frac{1}{i}\}$, and its Lyndon function coincides with the length $L$ of $%
\rho $-loops from $\Gamma (B)$, measured in $B$. Moreover, there is an
isometry of the $\mathbb{R}$-tree $T(G,L)$ onto $\overline{B}$, which
establishes the equivalence of the action of $G$ on $T(G,L)$ to the action
of $G=\Gamma (B)$ on $\overline{B}$.
\end{proposition}

\begin{proof}
Zastrow's group was defined by him as a subgroup of $\pi _{1}(H)\subset F$,
where $F$ is an inverse limit of a sequence $F_{n}$, $n\in \mathbb{N\,}$, of
free groups of rank $n$, using a complicated combinatorial description of $%
\pi _{1}(H)$ and $F$. But using the fact that by Proposition \ref{normal},
any element of $\pi _{1}(H)$ is represented by a normal loops based at $\ast 
$, we see that the first statement follows from the definitions of $G$ in
terms of the Lyndon function $L$ on $G$ (see p. 231, \cite{Ch}), $B$, and $%
\Gamma (B)$. Also $L=L_{\ast =1}$, by our definition, $\Gamma (B)\subset 
\overline{B}$, and $\Gamma (B)$ is the orbit of the point $\ast $ relative
to the isometric action of $\Gamma (B)$ on $\overline{B}$. It is clear that
any $\rho $-path $c$ in $B$ starting at $\ast $ is an initial part of a $%
\rho $-loop in $B$. This implies that the $\mathbb{R}$-tree $\overline{X}$
is \textit{minimal} in the above sense. By the above discussion, the proof
is finished.
\end{proof}

\begin{theorem}
\label{incl} If $f: X\rightarrow Y$ is an injective map such that for any
rectifiable path $c$ in $X$, the path $f\circ c$ in $Y$ is rectifiable, and $%
f$ topologically embeds the image of $c$ into $Y$, then the natural induced
map $f_{\ast}: \Gamma(X)\rightarrow \Gamma(Y)$ is an injective homomorphism.
If $f$ is a bijection such that $f^{-1}$ has the same properties as $f$
above, then the groups $\Gamma(X)$ and $\Gamma(Y)$ are isomorphic. In
particular, the last statement is true for any bi-Lipshitz map $f:
X\rightarrow Y$.
\end{theorem}

This theorem is an immediate corollary of the definitions. It gives some
sufficient but most likely not necessary conditions for the next open
question:

\begin{question}
\label{isomo}When are the groups $\Gamma (X)$ and $\Gamma (Y)$ isomorphic,
or, more specifically, when is $f_{\ast }$ an isomorphism?
\end{question}

The following lemma is a corollary of Proposition \ref{group} and the
discussion prior to it.

\begin{lemma}
\label{loop} If $c$ is a nontrivial $\rho $-loop in $X$ starting at $\ast $,
then there is a unique maximal (by inclusion) $\rho $-path $\alpha $ in $X$
starting at $\ast $ so that for some non-trivial $\rho $-loop $\beta $ in $X$
starting at the end of the path $\alpha $, $c=\alpha \star \beta \star
\alpha ^{-1}$. Then $\beta $ is also unique. In this situation, for any
nonzero integer $n$, $c^{n}=\alpha \star \beta ^{n}\alpha ^{-1}$, and $%
L(c^{n})=2L(\alpha )+|n|L(\beta )$ if we consider $c^{n}$ as an element of $%
\Gamma (X)$.
\end{lemma}

\begin{proposition}
\label{pleni}Let any $\rho $-path in $X$ starting at $\ast $ be an initial
part of a $\rho $-loop in $X$ and suppose there is a topological embedding $%
f:B\rightarrow X$, where $B$ is the same as in Proposition \ref{Za}, such
that for any rectifiable path $c$ in $B$, the path $f\circ c$ in $X$ is
rectifiable. Then $\Gamma (X)$ is an infinitely generated locally free group
that is not free and not a free product of surface groups and abelian
groups, but acts freely on the $\mathbb{R}$-tree $\overline{X}$. Moreover,
the $\mathbb{R}$-tree $\overline{X}$ is a minimal invariant subtree with
respect to this action.
\end{proposition}

\begin{proof}
Lemma \ref{loop} implies that for $1\neq g\in \Gamma (X)$, there is a
natural number $N$ such that if $g=h^{n}$, where $n$ is an integer and $h\in
\Gamma (X)$, then $|n|\leq N$. The group $\Gamma (X)$ is locally free by
Theorem \ref{ut1}. These two statements mean that all the statements of
Lemma 5.3.1 in \cite{Ch} are true for the group $\Gamma (X)$. The special
conditions for the map $f$ and Theorem \ref{incl} imply that $\Gamma (X)$
contains a group isomorphic to Zastrow's group $G$. In Lemma 5.3.3 from \cite%
{Ch} it is proved that $G$ is not free. Then the Nielsen-Schreier theorem
implies that $\Gamma (X)$ is not free. In Lemma 5.3.4 from \cite{Ch}, which
requires only the statements in Lemmas 5.3.1 and 5.3.3, it is proved that $G$
is not a free product of surface groups and abelian groups. Applying the
same proof, we get that $\Gamma (X)$ is not a free product of surface groups
and abelian groups. We proved in Theorem \ref{ut1} that $\Gamma (X)$ acts
freely by isometries on the $\mathbb{R}$-tree $\overline{X}$. The as yet
unused assumption implies, as in the proof of Proposition \ref{Za}, the last
statement.
\end{proof}

\begin{proof}[Proof of Theorem \protect\ref{plenitude2}]
The second statement is proved in Proposition \ref{Za}. Suppose we are given
two length metrics $d_{s},d_{t}$ on $H$, defined by sequences $%
\{s_{i}\},\{t_{i}\},i\in \mathbb{N}$. Then there is an increasing integer
sequence $k(i)$ such that $t_{k(i)}\leq s_{i}$ for all $i\in \mathbb{N}$. We
can define a 1-Lipshitz map $f:(H,d_{s})\rightarrow (H,d_{t})$, which is
also a topological embedding, by the requirement that $f(\ast )=\ast $ and $%
f|{C_{i}}:(C_{i},d_{s})\rightarrow (C_{k(i)},d_{t})$ is a bijective $%
(t_{k(i)}/s_{i})$-Lipshitz map. By Theorem \ref{incl}, this map induces an
injective homomorphism $f_{\ast }:\Gamma (H,d_{s})\rightarrow \Gamma
(H,d_{t})$. This proves the first statement.
\end{proof}

\begin{proof}[Proof of Theorem \protect\ref{plenitude1}]
The proof of Theorem \ref{universal}, Theorem \ref{plenitude2}, and
Proposition \ref{Za} imply that all these spaces satisfy the conditions of
Proposition \ref{pleni}. An application of this proposition finishes the
proof.
\end{proof}

\begin{definition}
\label{lt} A length space $X$ is local $\mathbb{R}$-tree at a point $x\in X$
if there is a number $r>0$ so that the closed ball $B(x,r)$ is an $\mathbb{R}
$-tree. The space $X$ is said to be local $\mathbb{R}$-tree if it is local $%
\mathbb{R}$-tree at any of its point.
\end{definition}

An example of a local $\mathbb{R}$-tree is given by the first case of
Example \ref{Zas}. Any traditional graph with some compatible length metric
is also a local $\mathbb{R}$-tree. The following three propositions can be
easily deduced from definitions.

\begin{proposition}
\label{lt1} If $X$ is a length space which is a local $\mathbb{R}$-tree,
then $\Gamma (X)=\pi _{1}(X)$.
\end{proposition}

\begin{proposition}
\label{analog} Let $(X,\ast )$ be the wedge product of length spaces $%
(X_{1},\ast )$ and $(X_{2},\ast )$ supplied with the natural length metric.
Then the group $\Gamma (X,\ast )$ is isomorphic to the free product $\Gamma
(X_{1},\ast )\ast \Gamma (X_{2},\ast )$ if at least one of $X_{1}$ or $X_{2}$
is local $\mathbb{R}$-tree at $\ast $.
\end{proposition}

\begin{proposition}
\label{wedge} Let $(X,d)$ be a length (respectively, geodesic) space which
is a local $\mathbb{R}$-tree at a point $\ast \in X.$ Then the closure in $X$
of any connected component of $X-\{\ast \}$ with the subspace metric $d$ is
a length (respectively, geodesic) space, $(X,\ast )$ is the wedge product of
the family $\{(X_{\alpha },\ast ),\alpha \in A\}$ of all such closures, and
the group $\Gamma (X,\ast )$ is isomorphic to the free product $\overset{%
\ast }{\prod }_{\alpha \in A}\Gamma (X_{\alpha },\ast )$ of the groups $%
\Gamma (X_{\alpha },\ast )$, $\alpha \in A$.
\end{proposition}

\begin{proposition}
\label{wedge1} For any family $\{X_{\alpha },\alpha \in A\}$ of length
(respectively, geodesic) spaces there exists a length (respectively,
geodesic) space $X$ such that the group $\Gamma (X)$ is isomorphic to the
free product $\overset{\ast }{\prod }_{\alpha \in A}\Gamma (X_{\alpha })$ of
the groups $\Gamma (X_{\alpha }),\alpha \in A$.
\end{proposition}

\begin{proof}
For every $\alpha \in A$, choose an arbitrary point $\star \in X_{\alpha }$
and let $(X_{\alpha }^{\prime },\ast )$ be $X_{\alpha }$ together with a
segment $\sigma _{\alpha }$ of fixed length $l>0$ with endpoints $\star $
and $\ast $, attached to $X_{\alpha }$ in such a way that $\sigma _{\alpha }$
has the point $\star $ in common with $X_{\alpha }$. By Proposition \ref%
{wedge}, $\Gamma (X_{\alpha }^{\prime },\star )$ is isomorphic to the free
product $\Gamma (X_{\alpha },\star )\ast \Gamma (\sigma _{\alpha },\star
)=\Gamma (X_{\alpha },\star )$ because $\Gamma (\sigma _{\alpha },\star )$
is the trivial group. By Proposition \ref{equivo}, the groups $\Gamma
(X_{\alpha }^{\prime },\star )$ and $\Gamma (X_{\alpha }^{\prime },\ast )$
are isomorphic. Then define $(X,\ast )$ as the wedge product of the spaces $%
(X_{\alpha }^{\prime },\ast )$. It is clear that $(X,\ast )$ is a length
(respectively, geodesic) space if all $X_{\alpha }$ are length
(respectively, geodesic) spaces, and $(X,\ast )$ is local $\mathbb{R}$-tree
at $\ast $. Furthermore, the closures of connected components of $(X,\ast
)-\{\ast \}$ are exactly the spaces $X_{\alpha }^{\prime }$, $\alpha \in A,$
and we can apply Proposition \ref{wedge}.
\end{proof}

\section{Piecewise continuously differentiable paths}

Let $(X,d)$ be any (connected) Riemannian manifold $M^{n}$ of dimension $%
n\geq 2$ with its length metric. Let $\tilde{X}$ consist of all $\rho $%
-paths in $(X,d)$ starting at $\ast \in X$ that are piecewise continuously
differentiable, and $\tilde{\Gamma}(X)\subset \tilde{X}$ be the
corresponding group of loops at $\ast $. Denote by $\star $ the restriction
of the operation $\star $ to $\tilde{\Gamma}(X)$.

\begin{theorem}
\label{utr}The $\mathbb{R}$-tree $\tilde{X}$ is a subtree of the $\mathbb{R}$%
-tree $(\overline{X},\overline{d})$ with induced metric and $(\tilde{\Gamma}%
(X),\star )$ is a (locally free) subgroup of $(\Gamma (X),\star )$. The $%
\mathbb{R}$-tree $\tilde{X}$ has the valency $\mathfrak{c}$ at each point
but is never complete. The length space $(X,d)$ is the metric quotient of $(%
\tilde{X},\tilde{d})$ via the free isometric action of the group $\tilde{%
\Gamma}(X)$ on $\tilde{X}$. The quotient mapping $\tilde{\phi}:\tilde{X}%
\rightarrow X$ is an arcwise isometry, a weak submetry (hence open) and
light map, and $\tilde{\phi}$ is a submetry if $X$ is geodesic. Moreover, $%
\tilde{X}$ is the minimal invariant subtree relative to the action of $%
\tilde{\Gamma}(X)$.
\end{theorem}

\begin{proof}
The first statement is evident, and this implies that the group $(\tilde{%
\Gamma}(X),\star )$ acts freely via isometries on $\tilde{X}$. Since $\tilde{%
X}$ is a subtree of the $\mathbb{R}$-tree $\overline{X}$, which is isometric
to $A_{\mathfrak{c}}$ by Theorem \ref{universal}, $\overline{X}$ has valency
at each point no more than $\mathfrak{c}.$ On the other hand, we can extend
any path $c\in \tilde{X}$ with the endpoint $c(L)$ by a geodesic segment
starting at $x=c(L)$ which has arbitrary tangent unit vector $v$ at the
point $x$. Since we have $\mathfrak{c}$ such vectors, Proposition \ref{val}
implies that $\overline{X}$ has valency $\mathfrak{c}$ at each point. It is
easy to construct a rectifiable map $c:[0,L)\rightarrow X$ starting at $\ast 
$ such that for any number $L_{k}$, $0<L_{k}<L$, the restriction $%
c_{k}=c|[0,L_{k}]$ is an arc-length parameterized piece-wise continuously
differentiable path, but $c$ either cannot extend continuously to some $%
c^{\prime }:[O,L]\rightarrow X$ (if $X$ is not complete), or such an
extension $c^{\prime }$ exists but is not a piecewise continuously
differentiable path. If we assume now that $L_{k}\nearrow L$, then $c_{k}$
is a Cauchy sequence in $\tilde{X}$ which has no limit in $\tilde{X}$, and
so $\tilde{X}$ is not complete. As a corollary of Theorem \ref{ut1}, $%
\overline{\phi }$ is 1-Lipschitz. Then its restriction $\tilde{\phi}$ is
also 1-Lipschitz. As in the proof of Theorem \ref{ut1}, using piecewise
continuously differentiable paths instead of more general rectifiable paths
in $X$, we get that the quotient mapping $\tilde{\phi}:\tilde{X}\rightarrow
X $ is an arcwise isometry, a weak submetry (hence open) and light map, and $%
\tilde{\phi}$ is a submetry if $X$ is geodesic. As in the proof of Theorem %
\ref{ut1}, this, together with the previously proved statements, implies the
third statement. The group $\tilde{\Gamma}(X)$ is locally free as a subgroup
of the locally free group $\Gamma (X)$. Considering only piecewise
continuous $\rho $-paths in $X$, we get from Proposition \ref{pleni} that $%
\tilde{X}$ is the minimal invariant subtree relative to the action of $%
\tilde{\Gamma}(X).$
\end{proof}

\begin{remark}
Notice that for $\tilde{\Gamma}(X)$ in the above theorem we may take also
the subset of $\overline{X}$, consisting of broken geodesics in $X$. Since $%
\tilde{\Gamma}(X)$ is locally free and satisfies condition (1) from Lemma
5.3.1 in \cite{Ch}, it cannot include a subgroup that is isomorphic to the
fundamental group of a surface or a non-cyclic subgroup of $(\mathbb{R},+)$.
So Question \ref{Mo1} for the group $\tilde{\Gamma}(X)$ is equivalent to the
question of whether $\tilde{\Gamma}(X)$ is a free group. We don't have any
answer to this question.
\end{remark}


\begin{thebibliography}{99}
\bibitem{AP} P.S. Alexandroff, B.A.Pasynkov, \textit{Introduction to
Dimension Theory} (Russian), Moscow, 1973.

\bibitem{AM} R.C. Alperin and K.N. Moss, Complete trees for groups with a
real-valued length function, J. London Math. Soc. (2) 31 (1985) 55-68.

\bibitem{An} R. D. Anderson, A continuous curve admitting monotone open maps
onto all locally connected metric continua, BAMS 62 (1956) 264-265.

\bibitem{An2} R. D. Anderson, Open mappings of continua, Summer Institute on
Set Theoretic Topology, Amer. Math. Soc., Providence, R. I., 1958.

\bibitem{AB} P. D. Andreev and V. N. Berestovski\v{i}, Dimensions of $%
\mathbb{R}$-trees and non-positively curved self-similar fractal spaces,
(Russian). Mat. Trudy 9 (2006), No.~2, 3-22.\newline
Engl.transl.: Siberian Adv. Math., 2007, v. 17, No.~2, 1-12.

\bibitem{Bs} V.N. Berestovski\v{i}, Homogeneous spaces with intrinsic
metric, Soviet Math. Dokl. 27 (1989) 60-63.

\bibitem{BG} V.N. Berestovski\v{i} and L. Guijarro, A Metric
Characterization of Riemannian Submersions, Annals of Global Analysis and
Geometry 18 (2000) 577-588.

\bibitem{BPU} V. Berestovski\v{i} and C. Plaut, Uniform universal covers of
uniform spaces, Top. Appl. 154 (2007) 1748-1777.

\bibitem{Be} M. Bestvina, $\mathbb{R}$-trees in topology, geometry and group
theory \textit{Handbook of Geometric Topology}, edited by R. Daverman and R.
Sher, Elsevier, Amsterdam, 2002, 55-91.

\bibitem{BF} M. Bestvina and M. Feighn, Stable actions of groups on real
trees, Inventiones Math. 121 (1995) 287-321.

\bibitem{B} R. H. Bing, Partitioning a set, Bull. Amer. Math. Soc. 55 (1949)
1101-1110.

\bibitem{BH} M. Bridson and A. Haefliger, \textit{Metric Spaces of
Non-positive Curvature}, Grundlehren der mathematishen Wissenshaften 319,
Springer-Verlag, Berlin, Heidelberg, New York, 1999.

\bibitem{Ch} I. Chiswell, \textit{Introduction to }$\Lambda $\textit{-trees}%
, World Scientific, Singapore, 2001.

\bibitem{CF} M. Curtis and M. Fort, Jr., The fundamental group of
one-dimensional spaces, Proc. Amer. Math. Soc. 10 (1959) 140--148.

\bibitem{CF1} M. Curtis and M. Fort, Jr., Homology of one-dimensional
spaces, Ann. Math., 2-nd ser. 69(2) (1959) 309--313.

\bibitem{Dun} M.J. Dunwoody, Groups acting on protrees, J. London Math. Soc.
56 (1997) 125-136.

\bibitem{DP} A. Dyubina and I. Polterovich, Explicit constructions of
universal $\mathbb{R}$-trees and asymptotic geometry of hyperbolic spaces,
Bull. Lon. Math. Soc. 33 (2001) 727-734.

\bibitem{F} K. J. Falconer, \textit{The Geometry of Fractal Sets}, Cambridge
University Press, Cambridge, 1985.

\bibitem{GLP} D. Gaboriau, G. Levitt and F. Paulin, Pseudogroups of
isometries of $\mathbb{R}$ and Rips' Theorem on free actions on $\mathbb{R}$%
-trees, Israel J. Math. 87 (1994) 403-428.

\bibitem{Gr} M. Gromov, \textit{Metric Structures for Riemannian and
Non-Riemannian Spaces}, Birkhauser, Boston, Basel, Berlin, 1999.

\bibitem{Hu} Bruce Hughes, Trees and ultrametric spaces: a categorial
equivalence, Adv. in Mathematics, 189(1) (2004) 148-191.

\bibitem{HW} W. Hurewicz, and H. Wallman, \textit{Dimension theory},
Princeton, 1948.

\bibitem{Ke} L. V. Keldy\v{s}, Example of a one-dimensional continuum with a
zero-dimensional and interior mapping onto a square (Russian), Doklady Akad.
Nauk SSSR (N.S.) 97 (1954) 201-204.

\bibitem{K} A. N. Kolmogorov, \"{U}ber offene abbildungen, Ann. of Math. (2)
38 (1937) 36--38.

\bibitem{KW} L. Kramer and R.M. Weiss, Coarse rigidity of Euclidean
buildings, preprint arXiv:0902.1332 [math.MG] 8 Feb 2009.

\bibitem{Ku} K. Kuratowski, \textit{Topology}, Academic Press. New York and
London, 1968.

\bibitem{MNO} J. Mayer, J. Nikiel, and L. Oversteegen, Universal spaces for $%
\mathbb{R}$-trees, TAMS 334 (1992) 411-432.

\bibitem{M} K. Menger, Untersuchungen \"uber allgemeine Metrik, Math. Ann.
vol. 100 (1928) 75-163.

\bibitem{Moi} E. Moise, Grille decomposition and convexification theorems
for compact locally connected continua, Bull. Amer. Math. Soc. (1949)
1111-1121.

\bibitem{Mo} J. W. Morgan, Deformations of algebraic and geometric
structures, CBMS Lectures, UCLA, summer 1986.

\bibitem{Mo1} J. W. Morgan, $\Lambda$-trees and their applications, Bull.
Amer. Math. Soc. (N.S.) 26(1) (1992) 87-112.

\bibitem{MS} J. W. Morgan and P. Shalen, Valuations, trees, and
degenerations of hyperbolic structures, I, Ann. Math 120 (1984) 401-476.

\bibitem{MS1} J. W. Morgan and P. Shalen, Free actions of surface groups on
R-trees, Topology 30 (1991), 143-154.

\bibitem{Po} C. Plaut, Metric spaces of curvature $\geq k$, Chapter 16, 
\textit{Handbook of Geometric Topology}, edited by R. Daverman and R. Sher,
Elsevier, Amsterdam, 2002.

\bibitem{T} J. Tits, A \textquotedblleft theorem of
Lie-Kolchin\textquotedblright\ for trees, in \textit{Contributions to
Algebra: A Collection of Papers Dedicated to Ellis Kolchin}, edited by H.
Bass, P. Cassidy, and J. Kovacic, Academic Press, New York, 1977.

\bibitem{W} D. Wilson, Open mappings of the universal curve onto continuous
curves, TAMS 168 (1972) 497-515.

\bibitem{Za} A. Zastrow, Construction of an infinitely generated group that
is not a free product of surface groups and abelian groups, but which acts
freely on an $\mathbb{R}$-tree, Proc. Royal Soc. Edinburg (A) 128(1998),
433-445.
\end{thebibliography}
\end{document}